\documentclass[reqno,12pt]{amsart}
\usepackage{amsfonts,amsmath,amsthm,amssymb,amscd}
\usepackage[all]{xy}
\newcommand{\version}{Ver.~0.0}
\newcommand{\setversion}[1]{\renewcommand{\version}{Ver.~{#1}}}
\setversion{1.0 [Wed Jan 24 15:22:16 JST 2007]}
\setversion{1.1 [Thu Jan 25 17:03:34 JST 2007]}

\title [Resolution and Conormal bundle]
{Resolution of null fiber and conormal bundles on the Lagrangian Grassmannian}

\dedicatory{Dedicated to Professor Noriaki Kawanaka on his 60th anniversary}

\author{Kyo Nishiyama}
\address{
Department of Mathematics\\
Graduate School of Science\\
Kyoto University\\
Sakyo, Kyoto 606-8502, Japan}
\email{kyo@math.kyoto-u.ac.jp}

\thanks{Partially supported by JSPS Grant-in-Aid for Scientific Research (B) \#{17340037}.}

\subjclass[2000]{Primary 14M15, 14E15; Secondary 14L30, 22E47, 20G05}
\keywords{resolution, null fiber, conormal bundle, Lagrangian Grassmannian, flag variety, dual pair, nilpotent orbit}

\theoremstyle{plain}
\newtheorem{theorem}{Theorem}
\newtheorem{proposition}[theorem]{Proposition}

\newtheorem{lemma}[theorem]{Lemma}

\newtheorem{introtheorem}{Theorem}

\theoremstyle{definition}

\theoremstyle{remark}
\newtheorem{remark}[theorem]{\upshape Remark}

\numberwithin{equation}{section}
\numberwithin{theorem}{section}

\newcommand{\Z}{\mathbb{Z}}

\newcommand{\R}{\mathbb{R}}

\newcommand{\bbK}{\mathbb{K}}
\newcommand{\bbL}{\mathbb{L}}

\newcommand{\bbG}{\mathbb{G}}
\newcommand{\C}{\mathbb{C}}

\newcommand{\bbU}{\mathbb{U}}
\newcommand{\bbV}{\mathbb{V}}

\newcommand{\lie}[1]{\mathfrak{#1}}

\newcommand{\Lie}{\mathop\mathrm{Lie}\nolimits{}}

\newcounter{thmenum}
\newenvironment{thmenumerate}{%
\begin{list}{$(\thethmenum)$}{%
\usecounter{thmenum}
\setlength{\labelsep}{.5em}
\setlength{\labelwidth}{-7pt}
\setlength{\topsep}{0pt}
\setlength{\partopsep}{0pt}
\setlength{\parsep}{0pt}
\setlength{\leftmargin}{3pt}
\setlength{\rightmargin}{0pt}
\setlength{\itemindent}{\leftmargin}
\setlength{\itemsep}{0pt}
}}
{\end{list}}

\newcommand{\mycomment}[1]{} 

\newlength{\lengthcup}
\newcommand{\rank}{\qopname\relax o{rank}}

\newcommand{\Hom}{\qopname\relax o{Hom}}
\newcommand{\Homfr}{\qopname\relax o{Hom}^{\circ}}

\newcommand{\Iso}{\qopname\relax o{Iso}}
\newcommand{\Stab}{\qopname\relax o{Stab}}

\newcommand{\id}{\qopname\relax o{id}}

\renewcommand{\Re}{\qopname\relax o{Re}}
\renewcommand{\Im}{\qopname\relax o{Im}}

\newcommand{\Image}{\qopname\relax o{Im}}

\newcommand{\Ker}{\qopname\relax o{Ker}}

\newcommand{\transpose}[1]{\,{}^t{#1}}

\newcommand{\closure}[1]{\overline{#1}}

\newcommand{\restrict}{\big|}
\newcommand{\partition}{\mathcal{P}}

\newcommand{\Spec}{\mathop\mathrm{Spec}\nolimits{}}
\newcommand{\Sym}{\mathop\mathrm{Sym}\nolimits{}}
\newcommand{\Alt}{\mathop\mathrm{Alt}\nolimits{}}

\newcommand{\orbit}{\mathbb{O}}  %

\newcommand{\length}[1]{\ell({#1})}

\newcommand{\nullcone}{\lie{N}}
\newcommand{\nullconeo}{\nullcone^{\circ}}
\newcommand{\resolution}{\widetilde{\lie{N}}}
\newcommand{\simresolution}{\widehat{\lie{N}}}

\newcommand{\nullconeofG}[1]{\nullcone_{\KGL}(#1)}

\newcommand{\harmonics}{\mathcal{H}}

\newcommand{\composit}{\odot}

\newcommand{\GITquotient}{/\!/}
\newcommand{\git}{\GITquotient}

\newcommand{\GL}{\mathrm{GL}}

\newcommand{\OO}{\mathrm{O}}
\newcommand{\U}{\mathrm{U}}

\newcommand{\Sp}{\mathrm{Sp}}
\newcommand{\USp}{\mathrm{USp}}
\newcommand{\Mat}{\mathrm{M}}

\newcommand{\Det}{\qopname\relax o{Det}}
\newcommand{\Grass}{\qopname\relax o{Grass}}

\newcommand{\fvariety}{\mathfrak{F}}
\newcommand{\Tbundle}{\mathcal{T}}
\newcommand{\Qbundle}{\mathcal{Q}}

\newcommand{\surjection}{\twoheadrightarrow}
\newcommand{\injection}{\hookrightarrow}

\newcommand{\LLp}{\langle\!\langle}
\newcommand{\LLm}{\langle\!\langle}
\newcommand{\LLpm}{\langle\!\langle}
\newcommand{\RRpm}{\rangle\!\rangle_{\scriptscriptstyle{\pm}}}
\newcommand{\RRp}{\rangle\!\rangle_{\scriptscriptstyle{+}}}
\newcommand{\RRm}{\rangle\!\rangle_{\scriptscriptstyle{-}}}

\newcommand{\XX}{\mathfrak{X}}
\newcommand{\XXp}{\mathfrak{X}^+}
\newcommand{\XXm}{\mathfrak{X}^-}

\newcommand{\orbitmin}{\orbit_{\mathrm{min}}}
\newcommand{\orbitminC}{\orbit_{\mathrm{min}\,\C}}

\newcommand{\KGL}{K}
\newcommand{\HGLxGL}{H}
\newcommand{\MOO}{K'}

\begin{document}

\begin{abstract} 
We study the null fiber of a moment map related to dual pairs.  
We construct an equivariant resolution of singularities of the null fiber, 
and get conormal bundles of closed $ K_{\C} $-orbits in the Lagrangian Grassmannian 
as the categorical quotient.  
The conormal bundles thus obtained turn out to be a resolution of singularities 
of the closure of nilpotent $ K_{\C} $-orbits, 
which is a ``quotient'' of the resolution of the null fiber.
\end{abstract}

\maketitle


\section*{Introduction}

Let $ \Mat_{m,n} = \Mat_{m,n}(\C) $ be the space of $ m \times n $-matrices with coefficients in $ \C $.  
Then the map $ \psi : \Mat_{m_1, n} \times \Mat_{n, m_2} \to \Mat_{m_1, m_2} $ defined by 
$ \psi(A, B) = A B $ is called a contraction map.  
There are many interpretations of this simple map.  
Let us explain one of them.

The general linear group $ \KGL = \GL_n(\C) $ acts on 
the space $ W = \Mat_{m_1, n} \times \Mat_{n, m_2} $ naturally by 
$ g \cdot ( A, B ) = ( A g, g^{-1} B ) \; ( g \in \KGL) $.  
Clearly the contraction map is invariant under this action of $ \KGL $.  
It turns out that $ \psi $ is the affine quotient map onto its image $ W \git K = \Spec( \C[W]^K ) $, 
which is the determinantal variety of rank $ r = \min \{ m_1, m_2, n \} $ (see Equation \eqref{eqn:determinantal.variety}).  

The fiber $ \nullconeofG{W} = \psi^{-1}(0) $ of the quotient map is called the \emph{null fiber} or \emph{null cone}.  
In the invariant theory, 
the study of the null fiber is important since it is the ``worst'' fiber of the quotient map.  
For example, if it has a finite number of orbits, then the other fibers also have a finite number of orbits.  
If it is reduced and irreducible, then so are the other fibers.  
For this, we refer the readers to \cite{Popov.Vinberg.1994}, for example.

In the following, we will give a natural construction of resolution of singularities for irreducible components of $ \nullconeofG{W} $, 
which is $ \KGL $-equivariant.  
There are many contributions to resolutions of null fibers, among them we refer  
\cite{Hesselink.1979} and \cite{Kraft.Wallach.2006} for general theory; and 
\cite{Buchsbaum.Eisenbud.1975}, \cite{Kempf.1975,Kempf.1976}, \cite{DeConcini.Strickland.1981}, \cite{Wang.2001}, etc.\ for 
the present situation (or more general treatment).  
Thus, the resolution itself has been well investigated.  

However, our main interest in this paper is 
to relate this resolution to the conormal bundle of a certain orbit on the Lagrangian flag variety associated to 
a symmetric pair.  
To state our main results precisely, we need to explain another interpretation of the contraction map.

Let us consider a dual pair $ ( \Sp_{2n}(\R), \OO_{m_1, m_2}(\R) ) $ in the large symplectic group 
$ \Sp_{2N}(\R) $ where $ N = n (m_1 + m_2) $.  
Then the vector space $ W $ above can be considered as a Lagrangian subspace in the symplectic space $ \C^N = \Mat_{2n, m_1 + m_2} $, 
and $ \psi $ is a moment map associated to the dual pair (see \S \ref{section:dual.pair}).  

Note that $ \KGL = \GL_n(\C) $ is the complexification of a maximal compact subgroup of $ \Sp_{2n}(\R) $, 
and if we take a Cartan decomposition of $ \lie{o}_{m_1, m_2}(\R) = \lie{k}'_{\R} \oplus \lie{p}'_{\R} $, 
the target space $ \Mat_{m_1, m_2} $ of the contraction map $ \psi $ 
is isomorphic to the complexification $ \lie{p}' = \lie{p}'_{\R} \otimes_{\R} \C $ of the Cartan space.

In completely analogous way, 
there also exists another moment map 
\begin{equation*}
\varphi : W \to \lie{p} = \Sym_n(\C) \oplus \Sym_n(\C)^{\ast} 
\qquad
\varphi(A, B) = (\transpose{A} A, B \transpose{B}) , 
\end{equation*}
where $ \lie{p} $ denotes (the complexification of) a Cartan space for $ \Sp_{2n}(\R) $.  
We can lift nilpotent orbits in $ \lie{p}' $ to those in $ \lie{p} $ using these moment maps 
(see Theorem \ref{theorem:NOZ} for precise statement).  
In this terminology, 
if we put $ \MOO = \OO_{m_1}(\C) \times \OO_{m_2}(\C) $ 
the complexification of a maximal compact subgroup of $ \OO_{m_1, m_2}(\R) $, 
the quotient of the null fiber $ \nullconeofG{W} $ by the action of $ \MOO $ gives a collection of 
the lifted spherical nilpotent orbits for the symmetric pair $ ( \Sp_{2n}(\R), \U_n(\C) ) $ which is 
usually denoted $ (+-)^p (-+)^q (+)^{n - (p + q)} (-)^{n - (p + q)} $ in terms of signed Young diagrams (cf.~\cite{Collingwood.McGovern.1993}).
Let us denote this orbit by $ \orbit_{p,q} $ (the precise definition will be given in \S \ref{section:dual.pair}).  
These two-step nilpotent orbits play an essential role in the following.

Now let us assume $ n = m_1 = m_2 $, i.e., we consider the \emph{equal rank case}.  
We keep this assumption till the end of Introduction.  
In equal rank case, the null fiber has the following irreducible decomposition
\begin{equation*}
\nullconeofG{W} = \textstyle\bigcup\nolimits_{p + q = n} \nullcone_{p,q} ,
\end{equation*}
with $ ( n + 1) $ irreducible components.  
Each component $ \nullcone_{p,q} \; ( p + q = n ) $ is stable under the action of $ \MOO $, and 
the categorical quotient by $ \MOO $ gives the closure of the nilpotent orbit $ \orbit_{p,q} $ mentioned above.  

We will construct a $ \KGL $-equivariant resolution of singularities of each component $ \nullcone_{p,q} $, 
which is denoted by $ \nu : \resolution_{p,q} \to \nullcone_{p,q} $.  
This resolution $ \resolution_{p,q} $ is realized as a vector bundle over a Grassmannian variety $ \Grass_p(\C^n) $ 
of $ p $-dimensional subspaces in $ \C^n $ (see \S \ref{subsection:outside.of.stable.range} and 
Proposition \ref{prop:resolution.outside.of.stable.range} for details).  
Let $ \XXm $ be the Lagrangian Grassmannian of all the Lagrangian subspaces in $ \C^{2n} $, 
which is isomorphic to a partial flag variety of $ \Sp_{2n}(\C) $ with respect to the Siegel parabolic subgroup.  
Then $ \KGL $ naturally acts on $ \XXm $, 
and it has exactly $ ( n + 1) $ closed orbits $ \mathcal{Z}_{p,q} \; ( p + q = n ) $ and 
$ \mathcal{Z}_{p,q} $ is isomorphic to $ \Grass_p(\C^n) $ (Lemma \ref{lemma:Zpq}).

With the setting explained above, the following theorems are main results of this article.

\begin{introtheorem}[Theorem \ref{theorem:conormal.bundle.is.quotient.of.reoslution}]
Assume $ n = m_1 = m_2 $, and let us denote $ \KGL = \GL_n(\C) $ and $ \MOO = \OO_{n}(\C) \times \OO_{n}(\C) $.
Then the categorical quotient $ \resolution_{p,q} \git \MOO $ of the resolution exists, and 
it is isomorphic to the conormal bundle $ T^{\ast}_{\mathcal{Z}_{p,q}} \XXm $ 
over the closed $ \KGL $-orbit $ \mathcal{Z}_{p,q} $ in the Lagrangian flag variety $ \XXm $.  
\end{introtheorem}

There is a natural moment map $ \mu $ from the cotangent bundle $ T^{\ast} \XXm $ to 
the nilpotent variety.  
If we restrict $ \mu $ to the conormal bundle $ T^{\ast}_{\mathcal{Z}_{p,q}} \XXm $ over $ \mathcal{Z}_{p,q} $, 
the image of the restricted $ \mu $ is the closure of a nilpotent $ \KGL $-orbit since there are only finite number of nilpotent orbits.

\begin{introtheorem}[Theorem \ref{theorem:resolution.of.nilpotent.orbit}]
The image of the moment map $ \mu : T_{\mathcal{Z}_{p,q}}^{\ast} \XXm \to \lie{p} $ coincides with 
the closure of the nilpotent orbit $ \closure{\orbit}_{p,q} \simeq \nullcone_{p,q} \git \MOO $.
Thus we have the following commutative diagram.
\begin{equation*}
\xymatrix{
\ar@{}[dr]|{\text{\Large$\circlearrowright$}}
\ar[d]_{\nu} \resolution_{p,q} \ar[r]^{\git \MOO} & T_{\mathcal{Z}_{p,q}}^{\ast} \XXm \ar[d]^{\mu} \\
\nullcone_{p,q} \ar[r]_{\git \MOO} & \closure{\orbit}_{p,q}
}
\end{equation*}
\end{introtheorem}

If $ n $ is even, there is an analogue of the above theorems for the dual pair 
$ ( \OO^{\ast}_{2n} , \Sp_{n.n} ) $.  
See \S \ref{section:O2n.ast.Spnn} for details.

\smallskip
The results of this article are purely geometric.  
However, we believe these geometrical properties of the null fiber reflect to the analysis of degenerate principal series  
of $ \Sp_{2n}(\R) $ and the theta liftings from $ \OO_{n,n}(\R) $ to $ \Sp_{2n}(\R) $.  

\medskip
\textbf{Acknowledgment.}  
The author thanks to Soo Teck Lee and Peter Trapa for useful discussions.

\section*{Notation}

Let $ \partition_n $ be the set of all the partitions $ \alpha = ( \alpha_1 , \alpha_2, \dots , \alpha_n ) , \;
\alpha_1 \geq \alpha_2 \geq \dots \geq \alpha_n \geq 0 $ with length less than or equal to $ n $.  
Here, the length of the partition $ \alpha $, denoted by $ \length{\alpha} $,  
is defined by the biggest number $ l $ such that 
$ \alpha_l > 0 = \alpha_{l + 1} $.  
We naturally identify $ \partition_n $ with the subset of $ \partition_{n + 1} $ 
by extending the sequence by adding $ 0 $.  

For $ \lambda \in \Z^n $, we say that $ \alpha $ is a dominant weight if 
$ \lambda_1 \geq \dots \geq \lambda_n $ holds.  
By definition, a partition $ \alpha \in \partition_n $ is always dominant.  
For two partitions $ \alpha \in \partition_s $ and 
$ \beta \in \partition_t $ with $ s + t \leq n $, 
we put 
\begin{equation*}
\alpha \composit \beta = ( \alpha , 0 , \dots , 0 , \beta^{\ast} ) 
= ( \alpha_1, \dots , \alpha_s , 0, \ldots, 0 , - \beta_t , \dots , - \beta_1 ) \in \Z^n ,
\end{equation*}
which is also a dominant weight.  
For a dominant $ \lambda $, we denote by $ \rho_{\lambda}^{(n)} $ 
the irreducible finite dimensional representation of $ \GL_n(\C) $ with the highest weight $ \lambda $.  
We often drop the superscript $ (n) $, if it is clear from the situation.  
The contragredient representation of $ \rho_{\lambda} $ is denoted by $ \rho_{\lambda}^{\ast} $, and 
the highest weight of $ \rho_{\lambda}^{\ast} $ by $ \lambda^{\ast} $.  

\section{Contraction map and its null fiber}
\label{section:contraction.map}

Let $ U_{i} = \C^{m_i} \; ( i = 1, 2 ) $ and $ V = \C^n $ be finite dimensional vector spaces, and put 
\begin{equation*}
W = V^{\ast} \otimes U_1 \oplus (U_2)^{\ast} \otimes V 
= \Hom (V, U_1) \times \Hom(U_2, V) .
\end{equation*}
We consider $ W $ as a representation space 
of $ \GL(U_1) \times \GL(U_2) \times \GL(V) $ in a natural manner.  
Put $ \HGLxGL = \GL(U_1) \times \GL(U_2) $ and $ \KGL = \GL(V) $ for short.  
The action of $ ( (h_1, h_2) , k ) \in \HGLxGL \times \KGL $ on $ ( f_1, f_2 ) \in \Hom (V, U_1) \times \Hom(U_2, V) $ is explicitly given by
\begin{equation*}
( (h_1, h_2) , k ) \cdot ( f_1, f_2 ) = ( h_1 \circ f_1 \circ k^{-1} , k \circ f_2 \circ h_2^{-1} )
\end{equation*}
The composite map 
\begin{equation}
\label{eqn:contraction.map}
\psi : W \to \Hom ( U_2, U_1 ) , \qquad \psi(f_1, f_2) = f_1 \circ f_2 
\end{equation}
is often called a contraction map.  
It is an affine quotient map by the action of $ \KGL $ onto its image.  
To state it more precisely, let $ r = \min \{ m_1, m_2, , n \} $.  
Then the image of $ \psi $ is 
\begin{equation}
\label{eqn:determinantal.variety}
\Det_r ( U_2 , U_1 ) = \{ f \in \Hom(U_2, U_1) \mid \rank f \leq r \} ,
\end{equation}
the determinantal variety of rank $ r $, so that $ W \git \KGL \simeq \Det_r(U_2, U_1) $.  

In general, suppose that a reductive algebraic group $ \KGL $ acts on an affine variety (or affine scheme) $ X = \Spec A $, 
where $ A $ is a commutative ring.  
Then an \emph{affine quotient} (or the categorical quotient in the category of affine schemes) is defined by 
$ X \git \KGL = \Spec A^{\KGL} $, where $ A^{\KGL} $ denotes the ring of $ \KGL $-invariants.  
Thus, in our case, an affine quotient is defined to be 
$ W \git \KGL = \Spec \C[W]^{\KGL} $ which is isomorphic to $ \Det_r(U_2, U_1) $ by the classical invariant theory.  
We abbreviate $ \Det_r(U_2, U_1) $ to $ \Det_r $, if there is no confusion.  

The fiber of $ 0 \in \Det_r $ is called the \emph{null fiber} or \emph{null cone}, and denoted as 
\begin{equation*}
\nullconeofG{W} = \{ ( f_1, f_2 ) \in W \mid f_1 \circ f_2 = 0 \}
= \{ ( f_1, f_2 ) \in W \mid \Image f_2 \subset \Ker f_1 \} .
\end{equation*}
It is easy to identify the $ \HGLxGL \times \KGL $-orbit structure of $ \nullconeofG{W} $.

\begin{lemma}
\label{lemma:orbit.decomposition.of.null.fiber}
Put $ \nullconeo_{s, t} = \{ ( f_1, f_2 ) \in \nullconeofG{W} \mid \rank f_1 = s , \rank f_2 = t \} $ 
and $ \nullcone_{s,t} = \closure{\nullconeo_{s,t}} $ its closure.  Then 
\begin{equation*}
\nullconeofG{W} = \textstyle\bigsqcup\nolimits_{
\begin{smallmatrix}
0 \leq s \leq m_1 , \, 0 \leq t \leq m_2 \\
 s + t \leq n 
\end{smallmatrix}
} \;
\nullconeo_{s,t}
\end{equation*}
gives the $ \HGLxGL \times \KGL $-orbit decomposition of the null fiber.  
We have 
\begin{equation*}
\dim \nullcone_{s, t} = \{ n - ( s + t ) \} ( s + t ) + s t + m_1 s + m_2 t ,
\end{equation*}
and the closure relation is described as follows; 
\begin{math}
\nullcone_{s', t'} \subset \nullcone_{s,t} 
\end{math}
if and only if $  s' \leq s \text{ and } t' \leq t $.  
Thus, the orbit $ \nullconeo_{p,q} $ is open in the null fiber if and only if $ p + q = \min \{ m_1 + m_2 , n \} $.
\end{lemma}

Let $ \harmonics_{\KGL} $ denote the space of \emph{harmonic polynomials} on $ W $, 
which are killed by the positive degree homogeneous $ \KGL $-invariant differential operators 
with constant coefficients (see, e.g., \cite{Howe.SchurLecture.1995} or \cite{Goodman.Wallach.1998} for precise definition).

\begin{lemma}
The space of harmonics decomposes as an $ \HGLxGL \times \KGL $-module in the following way.  
{\upshape(\emph{For the notation, we refer to those in the end of Introduction}.)}
\begin{equation*}
\harmonics_{\KGL} \simeq \textstyle \bigoplus_{
\begin{smallmatrix}
\alpha \in \partition_{m_1}, \, \beta \in \partition_{m_2} \\
\length{\alpha} + \length{\beta} \leq n
\end{smallmatrix}} \;
(\rho^{(m_1) \ast}_{\alpha} \boxtimes \rho^{(m_2)}_{\beta}) \boxtimes \rho^{(n)}_{\alpha \composit \beta}
\end{equation*}
\end{lemma}

Let $ J_{\KGL} $ be an ideal of the polynomial ring $ \C[W] $ 
generated by the positive degree homogeneous invariants $ \C[W]_+^{\KGL} $.  
$ J_{\KGL} $ is called an \emph{invariant ideal}.  
Then it is well known that 
\begin{equation*}
\C[W] = \harmonics_{\KGL} \oplus J_{\KGL} 
\qquad 
(\text{as an $ \HGLxGL \times \KGL $-module}) .
\end{equation*}
Since the ideal $ J_{\KGL} $ defines the null cone, 
there is a surjection from $ \harmonics_{\KGL} $ to the regular function ring over the null cone.  
Moreover, we have

\begin{theorem}
\label{theorem:decomposition.of.harmonics}
The null fiber $ \nullconeofG{W} $ is reduced, i.e., 
the invariant ideal $ J_{\KGL} $ is reduced.  
Thus we have $ \harmonics_{\KGL} \simeq \C[\nullconeofG{W}] $ as an $ \HGLxGL \times \KGL $-module.  
Moreover, we have
\begin{thmenumerate}
\item
If $ n \geq m_1 + m_2 $, the null fiber $ \nullconeofG{W} $ is irreducible, and 
it coincides with the closure of an open dense $ \HGLxGL \times \KGL $-orbit $ \nullconeo_{m_1,m_2} $.  
\item
If $ 1 \leq n < m_1 + m_2 $, the null fiber is reducible, 
and its irreducible decomposition is given by 
\begin{equation*}
\nullconeofG{W} = \textstyle\bigcup\nolimits_{
\begin{smallmatrix}
0 \leq p \leq m_1 , \, 0 \leq q \leq m_2 \\
 p + q = n 
\end{smallmatrix}
} \;\nullcone_{p,q} .
\end{equation*}
We have an irreducible decomposition of the $ \HGLxGL \times \KGL $-module 
\begin{equation*}
\C[ \nullcone_{p, q} ] \simeq 
\bigoplus_{\alpha \in \partition_p, \; \beta \in \partition_q}
(\rho^{(m_1) \ast}_{\alpha} \boxtimes \rho^{(m_2)}_{\beta}) \boxtimes \rho^{(n)}_{\alpha \composit \beta} 
\qquad
(\text{as an $ \HGLxGL \times \KGL $-module}).
\end{equation*}
\end{thmenumerate}
\end{theorem}

There are vast references concerning the null cone.  
The readers may refer \cite{Kraft.Wallach.2006} for the description of the irreducible component 
of the null cone in the general situation.  
The structure of the harmonics is essentially due to \cite{Kashiwara.Vergne.1978} and \cite{Howe.remarks.1989}.  
Also we refer to \cite[Theorem 3.3]{Brylinski.1993}.  

\section{Resolution of null fibers}

In this section, we briefly describe the construction of the resolution of each irreducible component 
$ \nullcone_{p,q} $ of the null fiber.  
The resolution has been well studied by many authors.  
We refer to the following articles 
and the references therein;   
\cite{Hesselink.1979}, \cite{Kraft.Wallach.2006} for general theory, 
and \cite{DeConcini.Strickland.1981}, \cite{Wang.2001}, etc.\ for the contraction map.   
So the content of this section is already known.  
However, because we need the precise construction and the notation later, we give the resolution explicitly in the following.

Let $ \nullconeo_{p,q} $ be an open orbit in the null fiber.  
Then it is equivariantly isomorphic to a fiber bundle over certain subvariety of a partial flag variety.  

\subsection{Stable range case ($ n \geq m_1 + m_2 $)}

First let us assume $ n \geq m_1 + m_2 $.  
In this case, the null fiber is irreducible and $ \nullconeo_{m_1,m_2} $ is 
the unique open $ \HGLxGL \times \KGL $-orbit in $ \nullconeofG{W} $.  
This open orbit has an explicit description, 
\begin{multline*}
\nullconeo_{m_1, m_2} = \{ ( f_1, f_2 ) \in \Hom(V, U_1) \times \Hom(U_2, V) \mid \\
\Image f_2 \subset \Ker f_1 , \; \rank f_1 = m_1 , \, \rank f_2 = m_2 \} ,
\end{multline*}
so that we can define a surjective map
\begin{equation*}
\zeta_{m_1, m_2} : \nullconeo_{m_1, m_2} \ni ( f_1 , f_2 ) \mapsto ( \Image f_2 , \Ker f_1 ) \in \fvariety_{(m_2, n - m_1)}(V) ,
\end{equation*}
where $ \fvariety_{(s,t)}(V) $ is a partial flag variety consisting of the pair of 
two subspaces $ ( V_2, V_1 ) $ of $ V $ 
which satisfies $ V_2 \subset V_1 $ and $ \dim V_2 = s , \; \dim V_1 = t $.  
We often abbreviate $ \fvariety_{(m_2, n - m_1)}(V) $ to $ \fvariety_{m_2, m_1'} $, where $ m_1' = n - m_1 $.  
The fiber $ \zeta_{m_1, m_2}^{-1}(x) $ at $ x = (V_2, V_1) \in \fvariety_{m_2, m_1'} $ is given by 
\begin{equation}
\label{eqn:fiber.of.zeta}
\zeta_{m_1, m_2}^{-1}(x) = \{ ( f_1, f_2 ) \mid 
f_1 : V / V_1 \xrightarrow{\sim} U_1 \text{ and } f_2 : U_2 \xrightarrow{\sim} V_2 \} , 
\end{equation}
where $ f_1 $ and $ f_2 $ are both isomorphisms, and $ f_1 $ is identified with the map 
$ V \surjection V / V_1 \xrightarrow{f_1} U_1 $, and similarly $ f_2 $ is identified with 
$ U_2 \xrightarrow{f_2} V_2 \injection V $.

Let us consider the tautological bundle $ \Tbundle_{m_2} $ on $ \fvariety_{m_2, m_1'} $.  
It is a vector subbundle of the trivial bundle $ \fvariety_{m_2, m_1'} \times V $ 
with the fiber $ V_2 $ at $ x = (V_2, V_1) \in \fvariety_{m_2, m_1'} $.  
Similarly, we can define $ \Tbundle_{m_1'} $ for which the fiber at $ x $ is $ V_1 $, 
and define the quotient bundle by 
\begin{equation*}
\Qbundle_{m_1'} = (\fvariety_{m_2, m_1'} \times V) / \Tbundle_{m_1'} .
\end{equation*}
Also we denote by $ \mathcal{U}_i $ the trivial bundle $ \fvariety_{m_2, m_1'} \times U_i $ ($ i = 1, 2 $).  
{From} the description of the fiber $ \zeta_{m_1, m_2}^{-1}(x) $ in  \eqref{eqn:fiber.of.zeta}, 
$ \nullconeo_{m_1, m_2} $ can be identified with the subbundle of the vector bundle 
$ \Hom(\Qbundle_{m_1'}, \mathcal{U}_1 ) \times \Hom(\mathcal{U}_2, \Tbundle_{m_2}) $.  
We denote this subbundle by 
$ \Iso(\Qbundle_{m_1'}, \mathcal{U}_1 ) \times \Iso(\mathcal{U}_2, \Tbundle_{m_2}) $.  
Namely it is explicitly given as
\begin{multline*}
\Iso(\Qbundle_{m_1'}, \mathcal{U}_1 ) \times \Iso(\mathcal{U}_2, \Tbundle_{m_2}) \\
= \{ ( x ; ( f_1, f_2 ) ) \in \Hom(\Qbundle_{m_1'}, \mathcal{U}_1 ) \times \Hom(\mathcal{U}_2, \Tbundle_{m_2}) \mid
\text{$ f_1 $ and $ f_2 $ are isomorphisms} \} .
\end{multline*}

\begin{proposition}
\label{prop:resolution.stable.range}
We assume the stable range condition $ n \geq m_1 + m_2 $, and put $ m_1' = n - m_1 $.  
\begin{thmenumerate}
\item
There is an $ \HGLxGL \times \KGL $-equivariant isomorphism
\begin{equation*}
\nullconeo_{m_1, m_2} \simeq \Iso(\Qbundle_{m_1'}, \mathcal{U}_1 ) \times \Iso(\mathcal{U}_2, \Tbundle_{m_2}) .
\end{equation*}
\item
The map 
\begin{equation*}
\nu_{m_1, m_2} : \Hom(\Qbundle_{m_1'}, \mathcal{U}_1 ) \times \Hom(\mathcal{U}_2, \Tbundle_{m_2}) \to \nullconeofG{W} , \quad
\nu_{m_1, m_2}( ( x ; (f_1, f_2 ) ) ) = (f_1 , f_2) 
\end{equation*}
gives an equivariant resolution of singularities of the null fiber.
\end{thmenumerate}
\end{proposition}

\begin{proof}
The first statement is obvious from the construction.  
Let us consider the second statement.  
We say the map $ \nu : X \to Y $ is a resolution of singularities of $ Y $ if 
(a) $ X $ is smooth, 
(b) $ \nu $ is a proper birational map, and 
(c) it is one to one precisely on the smooth points of $ Y $.  
In our case, (a) and the birationality in (b) are clear, 
and that the map $ \nu_{m_1, m_2} $ is proper 
follows because $ \Hom(\Qbundle_{m_1'}, \mathcal{U}_1 ) \times \Hom(\mathcal{U}_2, \Tbundle_{m_2}) $ 
is a vector bundle over a projective variety (see, e.g., \cite[Corollary 4.8 and Theorem 4.9]{Hartshorne.1977}).  
The condition (c) is also easy to verify.  
Note that $ \nullconeo_{m_1, m_2} $ is precisely the set of smooth points of the null fiber, 
which can be checked by the explicit calculations of the rank of differentials of defining equations.
\end{proof}

\subsection{Outside of the stable range ($ n < m_1 + m_2 $)}  
\label{subsection:outside.of.stable.range}

In this subsection, we assume that $ n < m_1 + m_2 $.  
Take an open $ \HGLxGL \times \KGL $-orbit $ \nullconeo_{p,q} \; ( p + q = n ) $ in the null fiber.  
In this case, for $ ( f_1, f_2 ) \in \nullconeo_{p,q} $, 
the image of $ f_2 $ and the kernel of $ f_1 $ coincide with each other since $ p + q = n $, 
and it gives a $ q $-dimensional subspace $ V_2 $ of $ V $.  
Thus, if we denote the Grassmann variety of $ q $-dimensional subspaces in $ V $ by $ \Grass_q(V) $, 
there is a $ \KGL $-equivariant map 
\begin{equation}
\zeta_{p, q} : \nullconeo_{p,q} \to \Grass_q(V) 
\end{equation}
which sends 
$ (f_1, f_2) \in \nullconeo_{p,q} $ to $ \Im f_2 = \Ker f_1 \in \Grass_q(V) $.  
This map is surjective, and $ \nullconeo_{p,q} $ is turned out to be a fiber bundle over $ \Grass_q(V) $ with 
the fiber at $ x = V_2 \in \Grass_q(V) $ 
\begin{equation*}
\zeta_{p,q}^{-1}(x) = \{ ( f_1, f_2 ) \mid 
f_1 : V / V_2 \injection U_1 \text{ and } f_2 : U_2  \surjection V_2 \} , 
\end{equation*}
where $ f_1 $ is injective and $ f_2 $ is surjective.  

We denote the tautological bundle on $ \Grass_q(V) $ by $ \Tbundle_q $, and the quotient bundle 
$ (\Grass_q(V) \times V) / \Tbundle_q $ by $ \Qbundle_q $.  As in the former subsection, 
$ \mathcal{U}_i = \Grass_q(V) \times U_i \; ( i = 1, 2 ) $ denotes the trivial bundle with the fiber $ U_i $.  
We denote by $ \Homfr( \Qbundle_q , \mathcal{U}_1 ) $ and $ \Homfr( \mathcal{U}_2 , \Tbundle_q ) $ 
the subbundle of the full rank maps in 
$ \Hom( \Qbundle_q , \mathcal{U}_1 ) $ and $ \Hom( \mathcal{U}_2 , \Tbundle_q ) $ respectively.

\begin{proposition}
\label{prop:resolution.outside.of.stable.range}
Let us assume that $ n < m_1 + m_2 $.  
\begin{thmenumerate}
\item
An open orbit $ \nullconeo_{p,q} \; ( p + q = n , \; 0 \leq p \leq m_1, \, 0 \leq q \leq m_2 ) $ in $ \nullconeofG{W} $ 
is $ \HGLxGL \times \KGL $-equivariantly isomorphic to 
$ \Homfr( \Qbundle_q , \mathcal{U}_1 ) \times \Homfr( \mathcal{U}_2 , \Tbundle_q ) $.
\item
The map 
\begin{equation*}
\nu_{p,q} : \Hom( \Qbundle_q , \mathcal{U}_1 ) \times \Hom( \mathcal{U}_2 , \Tbundle_q ) \to \nullcone_{p,q} , 
\quad
\nu_{p,q}( x ; ( f_1, f_2 ) ) = ( f_1, f_2 )
\end{equation*}
gives an $ \HGLxGL \times \KGL $-equivariant resolution of singularities of the irreducible component 
$ \nullcone_{p,q} $ of the null fiber.
\end{thmenumerate}
\end{proposition}

The proof of this proposition is completely similar to that of Proposition \ref{prop:resolution.stable.range}.  
We will denote the resolution of singularities in 
Propositions \ref{prop:resolution.stable.range} and \ref{prop:resolution.outside.of.stable.range} by 
\begin{equation}
\label{eqn:reolution.of.singularity.final.version}
\nu_{p,q} : \resolution_{p,q} \to \nullcone_{p,q} , \qquad
\resolution_{p,q} = \Hom( \Qbundle_{p'} , \mathcal{U}_1 ) \times \Hom( \mathcal{U}_2 , \Tbundle_q ) 
\end{equation}
in the following ($ p' = n - p $).

\subsection{Reduction to the stable range lift}
\label{subsection:reduction.to.stable.range}

Let us assume that we are still in the outside of the stable range, i.e., assume $ n < m_1 + m_2 $.  
The resolution of $ \nullcone_{p,q} \; ( p + q = n ) $ above can be reduced to the stable range case.  
Let us explain it.

For $ 0 \leq p \leq m_1, \; 0 \leq q \leq m_2 $ which satisfies $ n = p + q $, 
we take a subspace $ U^p \subset U_1, \; U^q \subset U_2 $ such that $ \dim U^p = p $ and $ \dim U^q = q $.  
Then one can play the same game for $ W^{p,q} = \Hom(V, U^p) \times \Hom(U^q , V ) $ instead of $ W $, 
with $ \HGLxGL = \GL(U_1) \times \GL(U_2) $ replaced by $ \GL(U^p) \times \GL(U^q) $.  
We can consider naturally $ W^{p,q} $ as a subspace of $ W $, 
and the intersection $ \nullcone_{p,q} \cap W^{p,q} $ is exactly the null fiber of the quotient map
$ W^{p,q} \to W^{p,q} \git \KGL = \Hom(U^q, U^p) $ in the stable range.  

Now we return to the $ \HGLxGL \times \KGL $-orbit $ \nullconeo_{p,q} $ in $ W $.  
To be sure, we again take $ \HGLxGL $ as $ \GL(U_1) \times \GL(U_2) $.  
For $ ( f_1, f_2 ) \in \nullconeo_{p,q} $, we can associate 
$ ( \Image f_1 , \Ker f_2 ) \in \Grass_p(U_1) \times \Grass_p(U_2) $.  
Notice that $ (\Ker f_2)^{\bot} = \Image f_2^{\ast} $, which is of dimension $ q = n - p $.  
The projection $ \nullconeo_{p,q} \to \Grass_p(U_1) \times \Grass_p(U_2) $ given above is surjective, 
and it gives $ \nullconeo_{p,q} $ a structure of a fiber bundle over the product of two Grassmannian.  
In this picture, $ \nullconeo_{p,q} $ is a subbundle of the vector bundle
\begin{equation*}
\Hom( \mathcal{V} , \Tbundle_p(U_1) ) \times \Hom( \Qbundle_p(U_2) , \mathcal{V} ) \to \Grass_p(U_1) \times \Grass_p(U_2) , 
\end{equation*}
where $ \mathcal{V} $ denotes the trivial bundle with the fiber $ V $.  
Thus $ \nullconeo_{p,q} $ can be identified with 
\begin{multline*}
\nullconeo_{p,q} = \{ ( x ; ( f_1, f_2 ) ) \in \Hom( \mathcal{V} , \Tbundle_p(U_1) ) \times \Hom( \Qbundle_p(U_2) , \mathcal{V} ) \mid \\
\text{$ f_1 $ and $ f_2 $ are of full rank and }
f_1 \circ f_2 = 0 \} ,
\end{multline*}
which means that it is a fiber bundle over $ \Grass_p(U_1) \times \Grass_p(U_2) $ with a fiber isomorphic to the dense orbit in the null fiber 
of the reduced quotient map $ W^{p,q} \to \Hom(U^q, U^p) $.  
We already know that the null fiber of the map $ W^{p,q} \to \Hom(U^q, U^p) $ has a resolution of singularities 
\begin{equation*}
\Hom( \Qbundle_{p'}(V) , U^p ) \times \Hom( U^q, \Tbundle_q(V) ) \to \Grass_q(V) 
\qquad (p' = n - p = q )
\end{equation*}
(see Equation~\eqref{eqn:reolution.of.singularity.final.version}).  
Thus we have a simultaneous resolution of $ \nullcone_{p,q} $ by the vector bundle
\begin{multline*}
\Hom( \Qbundle_{p'}(V) , \Tbundle_p(U_1) ) \times \Hom( \Qbundle_{q'}(U_2) , \Tbundle_q(V) ) \\
\to \Grass_q(V) \times\Grass_p(U_1) \times \Grass_p(U_2)  .
\end{multline*}
The explicit resolution map is given, for $ x = ( V_0 , V_1, V_2 ) \in \Grass_q(V)\times\Grass_p(U_1)\times\Grass_p(U_2) $, by
\begin{equation*}
( x ; ( F_1, F_2 ) ) \mapsto ( f_1, f_ 2 ) \quad \text{ where } \quad
\begin{cases}
f_1 : V \surjection V / V_0 \xrightarrow{F_1} V_1 \injection U_1  \\
f_2 : U_2 \surjection U_2 / V_2 \xrightarrow{F_2} V_0 \injection V
\end{cases}
\end{equation*}
We denote this simultaneous resolution by
\begin{equation}
\hat{\nu}_{p,q} : \simresolution_{p,q} 
= \Hom( \Qbundle_{p'}(V) , \Tbundle_p(U_1) ) \times \Hom( \Qbundle_{q'}(U_2) , \Tbundle_q(V) ) \to \nullcone_{p,q} . 
\end{equation}
Thus we have established the following theorem, which shows that indeed $ \simresolution_{p,q} $ is the collection of the resolutions of 
the null fibers in the stable range, which is parametrized by $ \Grass_p(U_1) \times \Grass_p(U_2) $.

\begin{theorem}
\begin{thmenumerate}
\item
There is a natural surjective, birational map $ f_{p,q} : \simresolution_{p,q} \to \resolution_{p,q} $ 
such that $ \hat{\nu}_{p,q} $ factors to 
$ \hat{\nu}_{p,q} = \nu_{p,q} \circ f_{p,q} $.
Moreover, the following diagram commutes.
\begin{equation*}
\begin{CD}
\simresolution_{p,q} @ > f_{p,q} > >  \resolution_{p,q}\\
@ V \hat{\nu}_{p,q} VV @ VV \nu_{p,q}V \\
\nullcone_{p,q} @ =  \nullcone_{p,q}
\end{CD}
\qquad\qquad
\begin{CD}
\simresolution_{p,q} @ > f_{p,q} > >  \resolution_{p,q}\\
@ V \tilde{\zeta}_{p,q} VV @ VV \zeta_{p,q}V \\
\Grass_{q,p}(V ; U_1, U_2 ) @ > \text{\upshape{projection}} > >  \Grass_p(V)
\end{CD}
\end{equation*}
where $ \Grass_{q,p}(V ; U_1, U_2 ) = \Grass_q(V) \times \Grass_p(U_1) \times \Grass_p(U_2) $.
\item
Fix a point $ x = ( X_1 , X_2 ) \in \Grass_p(U_1) \times \Grass_p(U_2) $ and put $ U^p = X_1 , \; U^q = U_2/X_2 $.  
Then 
\begin{equation*}
\nullcone_{p,q} \cap W^{p,q} = \{ ( f_1, f_2 ) \in \nullconeofG{W} \mid \Image f_1 \subset X_1 , X_2 \subset \Ker f_2 \}
\end{equation*}
is isomorphic to the null fiber of the contraction map $ W^{p,q} \to \Hom(U^p, U^q) $, 
and the map 
\begin{equation*}
\tilde{\zeta}_{p,q}^{-1} ( \Grass_q(V) \times \{ x \} ) \to \nullcone_{p,q} \cap W^{p,q} 
\end{equation*}
gives a resolution of singularities of the null cone $ \nullcone_{p,q} \cap W^{p,q} $.  
Moreover, $ f_{p,q} $ restricted to $ \tilde{\zeta}_{p,q}^{-1} ( \Grass_p(V) \times \{ x \} ) $ is an isomorphism 
into $ \resolution_{p,q} $.
\end{thmenumerate}
\end{theorem}

\section{Dual pair and moment maps}
\label{section:dual.pair}

In this section, we briefly recall the properties of moment maps and the theta lifting of nilpotent orbits.

Let us consider the dual pair $ ( G_{\R}, G_{\R}' ) = ( \Sp_{2n}(\R) , \OO_{m_1, m_2}(\R) ) $ in $ \bbG_{\R} = \Sp_{2N}(\R) $ where 
$ N = n ( m_1 + m_2 ) $ (see \cite{Howe.1985} for the definition of the dual pair).  
Here $ G_{\R} = \Sp_{2n}(\R) $ is realized as the group of linear transformations on $ \R^{2n} $ which preserve a given symplectic form.  
Also $ G_{\R}' = \OO_{m_1, m_2}(\R) $ is the indefinite orthogonal group realized on 
a quadratic space $ \R^{m_1, m_2} $ with the signature $ ( m_1, m_2 ) $.  
Then the tensor product $ W = \R^{2 n} \otimes \R^{m_1, m_2} $ has a natural symplectic form 
which is the product of the symplectic form on $ \R^{2n} $ and the symmetric form on $ \R^{m_1, m_2} $.  
We denote this symplectic form by $ \langle \, , \, \rangle $.  
Thus $ \bbG_{\R} = \Sp_{2N}(\R) $ defined by $ \langle \, , \, \rangle $ contains 
$ \Sp_{2n}(\R) $ and $ \OO_{m_1, m_2}(\R) $ as subgroups and they are mutually commutant to each other, i.e., they form a dual pair.

Let us introduce a complex structure on the symplectic space $ W $.  

Take a complete polarization $ W = X \oplus Y $.  
Since $ X $ and $ Y $ are totally isotropic and the symplectic form is non-degenerate, 
there is a complete pairing on $ X \times Y $.  
Thus we can identify $ Y = X^{\ast} $.  
Let us fix a standard Euclidean inner product on $ X $.  
Then it also gives an identification $ X \simeq X^{\ast} $, and we get an isomorphism 
$ X \to X^{\ast} \to Y $.  We denote this isomorphism by $ \eta $.  

Let us denote the complexification of $ W $ by $ W_{\C} = W \otimes_{\R} \C $ 
(we use the same notation for the complexification for any real vector space).  
Then $ \eta : X \to Y $ extends to $ \eta : X_{\C} \to Y_{\C} $ linearly.  
Put
\begin{equation}
L^+ = \{ x + i \eta(x) \mid x \in X_{\C} \} .
\end{equation}
It is easy to check that $ L^+ $ is Lagrangian in $ W_{\C} $ with respect to the complexified symplectic form.  
Then we get an isomorphism $ \gamma : L^+ \to W $ over $ \R $ defined by 
$ \gamma(z) = \Re z $, where $ \Re z $ denotes the real part of $ z $.  
By this isomorphism, we can identify $ X_{\C} \simeq L^+ \simeq W $, and
$ W $ enjoys a complex structure coming from this identification.  
Also we define a Hermitian inner product on $ X_{\C} $ using the Euclidean inner product on $ X $, 
and transfer it to $ W $.  
Denote this inner product by $ ( \, , \, ) $.  
Then the symplectic form  $ \langle \, , \, \rangle $ on $ W $ is recovered 
from it as $ \langle v, w \rangle = - \Im ( v, w ) $.  
Thus, the unitary group $ \U(W) $ with respect to $ ( \, , \, ) $ preserves the symplectic form, hence 
it is a maximal compact subgroup $ \bbK_{\R} $ of $ \bbG_{\R} $.  
Also its complexification $ \bbK \simeq \GL_N(\C) $ naturally acts on $ W $, $ X_{\C} $ and $ L^+ $ through the above isomorphism, 
and this in turn determines an embedding of $ \bbK $ into $ \bbG := \Sp_{2N}(\C) $.

Since the vector representation of $ \Sp_{2N}(\C) $ on $ W_{\C} $ is Hamiltonian, 
there is a moment map $ \mu : W_{\C} \to \lie{sp}_{2N}(\C)^{\ast} $, the algebraic dual of the Lie algebra of $ \Sp_{2N}(\C) $.  
The map $ \mu $ can be described explicitly by the formula 
$ \mu(w)(X) = \frac{1}{2} \langle X w , w \rangle \; ( w \in W , \, X \in \lie{sp}_{2N}(\C) ) $.  
Then it is easy to see that the image $ \Image \mu $ coincides with the closure of the minimal nilpotent coadjoint orbit 
$ \closure{\orbitminC} $.  
The moment map restricted to $ L^+ $ has an image which is orthogonal to $ \Lie (\bbK_{\R})_{\C} = \lie{K} $.  
If we denote the complexified Cartan decomposition with respect to $ \bbK_{\R} $ by $ \lie{G} = \lie{K} \oplus \lie{P} $, 
then 
\begin{equation}
\mu(L^+) = \closure{\orbitmin^+} \subset \closure{\orbitminC} \cap \lie{P}^{\ast} ,
\end{equation}
where $ \orbitmin^+ $ is one of the minimal nilpotent $ \bbK $-orbits in $ \lie{P}^{\ast} $.
Note that $ \mu $ is $ \bbK $-equivariant.

Now take a complete polarization $ \R^{2n} = X_0 \oplus Y_0 $.  
Then $ X = X_0 \otimes \R^{m_1, m_2} $ is a Lagrangian subspace in $ W $ and, 
if we define $ Y $ similarly, $ W = X \oplus Y $ gives a complete polarization.  
Using this $ X $, we follow the above recipe to construct a moment map.
We take maximal compact subgroups $ K_{\R}, K_{\R}' $ of $ G_{\R}, G_{\R}'$ respectively inside $ \bbK_{\R} $, 
hence we get a compatible (complexified) Cartan decomposition 
\begin{equation}
\lie{g} = \lie{k} \oplus \lie{p} , \quad
\lie{k} \subset \lie{K} , \quad 
\lie{p} \subset \lie{P} .
\end{equation}
Similarly we define $ \lie{g'} = \lie{k'} \oplus \lie{p'} $ to be a Cartan decomposition for $ G_{\R}' $.
Then we can define a map 
\begin{equation*}
\varphi : W \simeq X_{\C} \xrightarrow{\sim} L^+ \xrightarrow{\mu} \closure{\orbitmin^+} \xrightarrow{\text{res.}} \lie{p}^{\ast}
\end{equation*}
where the last map is the restriction of $ \lie{P}^{\ast} $ to the subspace $ \lie{p} $.  
This is a $ K := K_{\C} $-equivariant map.  
Since $ K' := K'_{\C} $ acts trivially on $ \lie{p} $, the map $ \varphi $ is $ K' $-invariant.  
Similarly we can define a $ K' $-equivariant map 
\begin{math}
\psi : W \to \lie{p'}^{\ast}
\end{math}
which is $ K $-invariant.  
Thus we have a double fibration of the complex vector space $ W $.

\begin{remark}
Since the action of $ G_{\R} $ or $ G_{\R}' $ on $ W $ is Hamiltonian, we can define moment maps 
$ \mu : W_{\C} \to \lie{g}^{\ast} $ (or $ \lie{g'}^{\ast} $ respectively).  
If we restrict this moment map to $ W \simeq L^+ $, then they coincide with $ \varphi $ or $ \psi $ respectively.
\end{remark}

If we denote $ V = (X_0)_{\C} = \C^n $ and $ U_i = \C^{m_i} \; (i = 1, 2) $ as in \S \ref{section:contraction.map}, 
we can identify $ W = \Hom(V, U_1) \oplus \Hom(U_2, V) $ and the map $ \psi $ coincides with 
the contraction map (see \eqref{eqn:contraction.map}).  

Let us describe the map $ \varphi : W \to \lie{p}^{\ast} $.  
We shall identify $ \lie{p}^{\ast} $ and $ \lie{p} $ by the Killing form, 
and note that there is a $ K $-stable decomposition $ \lie{p} \simeq \lie{p}_+ \oplus \lie{p}_- $, 
where $ \lie{p}_+ \simeq \Sym_n(\C) = \Sym(V) $ and $ \lie{p}_- \simeq \Sym(V^{\ast}) $.  
Here we denote by $ \Sym_n(\C) $ the space of $ n \times n $-symmetric matrices, 
and $ \Sym(V) $ denotes the space of symmetric tensors of degree two.
Take $ ( f_1, f_2 ) \in \Hom(V, U_1) \oplus \Hom(U_2, V) = W $.  
We identify $ U_i $ with $ U_i^{\ast} $ by the symmetric form on $ U_i \; (i = 1, 2) $, which comes from 
the complexification of the quadratic form on $ \R^{m_1, m_2} $.  Thus we have $ f_1^{\ast} \circ f_1 $ by composing 
\begin{equation*}
V \xrightarrow{ f_1 } U_1 \simeq U_1^{\ast} \xrightarrow{ f_1^{\ast} } V^{\ast} .
\end{equation*}
It is clear that $ f_1^{\ast} \circ f_1 \in \Sym(V^{\ast}) $.  Similarly we get $ f_2 \circ f_2^{\ast} \in \Sym(V) $.  
Now $ \varphi $ is given by 
\begin{equation}
\varphi( f_1 , f_2 ) = ( f_1^{\ast} \circ f_1 , f_2 \circ f_2^{\ast} ) \in \Sym(V^{\ast}) \oplus \Sym(V) .
\end{equation}
The map $ \varphi $ is an affine quotient map from $ W $ onto its image by the action of 
$ K' = \OO_{m_1}(\C) \times \OO_{m_2}(\C) $.  
The image is the determinantal variety in the symmetric matrices (or tensors).

\begin{theorem}
\label{theorem:NOZ}
Assume the stable range condition $ n \geq m_1 + m_2 $.  
For any $ K' $-orbit $ \orbit' $ in $ \lie{p}' $, 
there exists a single $ K $-orbit $ \orbit $ in $ \lie{p} $ such that 
$ \varphi(\psi^{-1}(\closure{\orbit'}) = \closure{\orbit} $ holds.  
This correspondence, called $ \theta $-lifting, establishes an injection from 
the set of $ K' $-orbits into that of $ K $-orbits: $ \theta : \lie{p}' / K' \injection \lie{p} / K $.  
Moreover, we have
\begin{thmenumerate}
\item
The $ \theta $-lifting carries a nilpotent $ K' $-orbit to a nilpotent $ K $-orbit.
\item
The affine quotient $ \psi^{-1}(\closure{\orbit'}) \git K' $ is isomorphic to $ \closure{\orbit} $.
\end{thmenumerate}
\end{theorem}

This theorem (in wider context) is obtained independently by \cite{NOZ.2006} and \cite{DKP.2002, DKP.2005} for nilpotent orbits.  
Also it is noticed by Takuya Ohta (private communication) and the author (\cite{Nishiyama.Tuebingen.2005}) that the lifting map 
can be defined for arbitrary orbits.  

For $ p + q \leq n $, we denote by $ \orbit_{p,q} $ the $ \theta $-lifting of the trivial orbit $ \{ 0 \} \subset \lie{p}' $.  
This is a two step nilpotent $ K $-orbit of the shape $ 2^{p + q} \cdot 1^{2 (n - ( p + q) )} $.  
The signed Young diagram is given as $ (+-)^p(-+)^q(+)^{n - (p + q)}(-)^{n - (p + q)} $.  
(See \cite{Nishiyama.Tuebingen.2000} for the more detailed properties of this orbit.)

\section{Intersection of Lagrangian flag varieties}

From now on, we assume $ n = p = q $, i.e., $ \dim V $ is equal to $ \dim U_1 = \dim U_2 $.  
This case is called \emph{equal rank case}.  
Henceforth, we will denote $ U = U_1 = U_2 $.

First we introduce two nondegenerate bilinear forms on $ \bbV := V \oplus V^{\ast} $, 
where $ V^{\ast} = \Hom( V, \C ) $ is the algebraic dual.  
One $ \LLp \, , \, \RRp $ of the two bilinear forms is a symmetric form, 
and the other $ \LLm \, , \, \RRm $ is a symplectic (or alternating) one.  
We put 
\begin{equation*}
\LLpm \xi , v \RRpm = \xi(v) \qquad
\text{ for $ \xi \in V^{\ast} $ and $ v \in V $.}
\end{equation*}
Also the subspaces $ V $ and $ V^{\ast} $ are assumed to be totally isotropic, so that 
$ \bbV = V \oplus V^{\ast} $ is a joint polar decomposition for the both $ \LLpm \, , \, \RRpm $.  
Recall that $ W = \Hom(V, U) \times \Hom(U, V) $.

\begin{lemma}
\label{lemma:null.fiber.iff.isotropic}
The pair of maps $ ( f_1 , f_2 ) \in W $ belongs to the null fiber $ \nullconeofG{W} $ if and only if 
$ \Image f_2 \oplus \Image f_1^{\ast} \subset \bbV $ is a joint isotropic subspace with respect to 
$ \LLpm \, , \, \RRpm $, where $ f_1^{\ast} : U^{\ast} \to V^{\ast} $ is the adjoint map of $ f_1 $.
\end{lemma}

For $ w = ( f_1, f_2 ) \in W $, we put $ \bbU_w = \Image f_2 \oplus \Image f_1^{\ast} $, 
which is a subspace of $ \bbV $.  
According to the lemma above, $ \bbU_w $ is totally isotropic for both $ \LLpm \, , \, \RRpm $ 
if and only if $ w \in \nullconeofG{W} $.  
Moreover, $ \bbU_w $ is a \emph{joint Lagrangian subspace} (i.e., a maximal totally isotropic space 
for both of $ \LLpm \, , \, \RRpm $) if and only if 
$ w $ belongs to one of the open orbit $ \nullconeo_{p,q} $ where $ p + q = n $.  

Let us denote by $ \mathcal{Z} $ a closed subvariety of $ \Grass_n(\bbV) $ of all the joint Lagrangian subspaces in $ \bbV $.  

\begin{lemma}
\label{lemma:Zpq}
\begin{thmenumerate}
\item
\label{lemma:Zpq:item:1}
If $ \bbL $ is a joint Lagrangian subspace of $ \bbV $, 
it is decomposed into a direct sum $ \bbL = \bbL \cap V \oplus \bbL \cap V^{\ast} $.
Thus we have a disjoint decomposition of $ \mathcal{Z} $ as 
\begin{equation*}
\mathcal{Z} = {\textstyle\bigsqcup_{p + q = n} \mathcal{Z}_{p,q}} , \qquad
\mathcal{Z}_{p,q} = \{ \bbL \in \mathcal{Z} \mid \dim \bbL \cap V = p , \; \dim \bbL \cap V^{\ast} = q \} .
\end{equation*}
Each $ \mathcal{Z}_{p,q} $ is a single closed $ \KGL $-orbit, hence an irreducible smooth projective variety.  
Moreover, there is a $ \KGL $-equivariant isomorphism 
\begin{equation*}
\mathcal{Z}_{p,q} \xrightarrow{\; \sim \;} \Grass_p(V) , \qquad
\bbL \mapsto \bbL \cap V .
\end{equation*}
\item
\label{lemma:Zpq:item:2}
The map $ z_{p,q} : \nullconeo_{p,q} \to \mathcal{Z}_{p,q} $ defined by $ z_{p,q}(w) = \bbU_w \; ( w \in \nullconeo_{p,q} ) $ is surjective.  
\end{thmenumerate}
\end{lemma}

\begin{proof}
Let $ \bbL $ be a joint Lagrangian subspace.  
We define an involution $ \tau $ on $ \bbV $ by $ \tau \restrict_V = \id_V $ and $ \tau \restrict_{V^{\ast}} = - \id_{V^{\ast}} $.  
If $ \bbL $ is stable under $ \tau $, then $ \bbL $ is homogeneous with respect to the direct sum decomposition $ \bbV = V \oplus V^{\ast} $, 
and we get the decomposition $ \bbL = \bbL \cap V \oplus \bbL \cap V^{\ast} $.  
Take $ u + \xi , v + \eta \in \bbL $, where $ u, v \in V $ and $ \xi, \eta \in V^{\ast} $.  
Then we have
\begin{align*}
\LLm u + \xi , \tau(v + \eta) \RRm 
&
= \LLm u + \xi , v - \eta \RRm 
= - \LLm u , \eta \RRm + \LLm \xi , v \RRm 
\\
&
= \eta(u) + \xi(v) 
= \LLp u, \eta \RRp + \LLp \xi , v \RRp
= \LLp u + \xi , v + \eta \RRp = 0
\end{align*}
Thus $ \tau(v + \eta) $ is orthogonal to $ \bbL $ with respect to the symplectic form $ \LLm \, , \, \RRm $.  
Since $ \bbL $ is maximally isotropic, this implies $ \tau(v + \eta) \in \bbL $.

Let us write $ \bbL_E = \bbL \cap E $ for a subspace $ E $ in $ \bbV $.  
Then an $ n $-dimensional subspace $ \bbL = \bbL_V \oplus \bbL_{V^{\ast}} $ is joint Lagrangian 
if and only if $ \bbL_{V^{\ast}} = ( \bbL_V )^{\bot} $.  
Thus $ \mathcal{Z}_{p,q} $ is isomorphic to $ \Grass_p(V) $, which proves it is a single $ \KGL $-orbit at the same time.  
Since $ \mathcal{Z}_{p,q} = \mathcal{Z} \cap \Grass_p(V) \times \Grass_q(V^{\ast}) $, it is closed.
Now we have proved \eqref{lemma:Zpq:item:1}.

The second assertion is clear from Lemma \ref{lemma:null.fiber.iff.isotropic}.
\end{proof}

Notice that $ \XX := \Grass_n(\bbV) \simeq \GL(\bbV) / P_n $, where $ P_n = \Stab_{\GL(\bbV)}(V) $ denotes a maximal parabolic subgroup which stabilizes $ V $.  
Then $ \XXp := \OO(\bbV) \cdot e P_n \subset \XX $ and $ \XXm := \Sp(\bbV) \cdot e P_n $ are also partial flag varieties 
consisting of all the Lagrangian subspaces.  
Thus there exist natural embeddings 
\begin{equation*}
\OO(\bbV) / P_n^+ \simeq \XXp \injection \XX , \qquad
\Sp(\bbV) / P_n^- \simeq \XXm \injection \XX , 
\end{equation*}
where $ P_n^+ = \Stab_{\OO(\bbV)}(V) $ and $ P_n^- = \Stab_{\Sp(\bbV)}(V) $ are maximal parabolic subgroups respectively.
It is clear $ \mathcal{Z} = \XXp \cap \XXm $ 
so that $ \mathcal{Z}_{p,q} \simeq \Grass_p(V) $ can be regarded as a closed $ \KGL $-orbit in both $ \XXp $ and $ \XXm $.
Note that
\begin{equation*}
\KGL = \GL(V) \simeq \OO(\bbV) \cap \Sp(\bbV) ,
\end{equation*}
and that 
$ ( \OO(\bbV), \KGL ) $ and $ ( \Sp(\bbV) , \KGL ) $ are complex symmetric pairs corresponding to 
real symmetric pairs $ ( \OO^{\ast}_{2n}, \U_n(\C) ) $ and $ (\Sp_{2n}(\R) , \U_n(\C) ) $ respectively.

\section{Case of $ \Sp_{2 n}(\R) \times \OO_{n, n}(\R) $}
\label{section:Sp.O}

We consider the standard symmetric bilinear form on $ U (= U_1) $ induced by the identity matrix.  
Then we can define the orthogonal group $ \OO(U) $ with respect to the symmetric form.  
By this non-degenerate symmetric form, we can identify $ U $ and $ U^{\ast} $, 
and we put $ U_2 = U^{\ast} $.  

As before we put $ \HGLxGL = \GL(U) \times \GL(U^{\ast}) $, and 
$ \MOO := \OO(U) \times \OO(U^{\ast}) \subset \HGLxGL $ 
is the complexification of a maximal compact subgroup $ \OO_n(\R) \times \OO_n(\R) $ of 
an indefinite orthogonal group $ \OO_{n,n}(\R) $.  
Similarly, in this section, 
$ \KGL = \GL_n(\C) $ should be regarded as the complexification of a maximal compact subgroup $ \U_n(\C) $ of $ \Sp_{2n}(\R) $, 
and we use the setting for the dual pair $ ( \Sp_{2 n}(\R) , \OO_{n, n}(\R) ) $.

Recall the $ \HGLxGL \times \KGL $-equivariant resolution of singularities 
\begin{equation}
\resolution_{p,q} = \Hom( \Qbundle_p , \mathcal{U} ) \times \Hom( \mathcal{U} , \Tbundle_p ) \to \nullcone_{p,q} 
\end{equation}
of the irreducible component $ \nullcone_{p,q} $ of the null fiber ($ p + q = n $).  
Since $ \HGLxGL $ acts on the vector bundle $ \Hom( \Qbundle_p , \mathcal{U} ) \times \Hom( \mathcal{U} , \Tbundle_p ) $ fiber wise, 
we can take a categorical quotient by the action of $ \MOO = \OO(U) \times \OO(U^{\ast}) $ easily, and get
\begin{equation*}
\Hom( \Qbundle_p , \mathcal{U} ) \times \Hom( \mathcal{U} , \Tbundle_p ) \git \MOO 
\simeq 
\Sym( \Qbundle_p^{\ast} ) \times \Sym( \Tbundle_p ) ,
\end{equation*}
where $ \Sym(X) $ denotes the symmetric tensor of $ X $ in $ X \otimes X $.  
It is a vector bundle over $ \Grass_p(V) \simeq \mathcal{Z}_{p,q} $ on which $ \KGL = \GL(V) $ acts.

\begin{theorem}
\label{theorem:conormal.bundle.is.quotient.of.reoslution}
Let us denote the conormal bundle of $ \mathcal{Z}_{p,q} \subset \XXm $ by $ T_{\mathcal{Z}_{p,q}}^{\ast} \XXm $.  
Then there exists a $ \KGL $-equivariant isomorphism 
\begin{equation*}
T_{\mathcal{Z}_{p,q}}^{\ast} \XXm \simeq \Sym( \Qbundle_p^{\ast} ) \times \Sym( \Tbundle_p ) \simeq \resolution_{p,q} \GITquotient \MOO  .
\end{equation*}
Thus the conormal bundle of the closed $ \KGL $-orbit $ \mathcal{Z}_{p,q} $ in $ \XXm $ is obtained by 
taking the categorical quotient of the resolution of singularities $ \resolution_{p,q} $ of $ \nullcone_{p,q} $.
\end{theorem}

\begin{proof}
We have to prove $ T_{\mathcal{Z}_{p,q}}^{\ast} \XXm \simeq \Sym( \Qbundle_p^{\ast} ) \times \Sym( \Tbundle_p ) $.  
Let us choose a point $ \bbL \in \mathcal{Z}_{p,q} $.  Then the tangent space of $ \XXm $ at $ \bbL $ can be identified with 
$ \Sym(\bbL^{\ast}) $, 
while the tangent space of $ \mathcal{Z}_{p,q} $ at $ \bbL = \bbL_V \oplus \bbL_{V^{\ast}} $ should be 
\begin{equation*}
\Hom(\bbL_V, V/\bbL_V) \simeq \bbL_V^{\ast} \otimes (V/\bbL_{V}) \simeq \bbL_V^{\ast} \otimes \bbL_{V^{\ast}}^{\ast} 
\end{equation*}
(recall that $ \bbL_{V^{\ast}} = (\bbL_V)^{\bot} $).  
Since 
$ \Sym(\bbL^{\ast}) \simeq \Sym(\bbL_V^{\ast}) \oplus \Sym(\bbL_{V^{\ast}}^{\ast}) \oplus \bbL_V^{\ast} \otimes \bbL_{V^{\ast}}^{\ast} $, 
we have 
\begin{equation*}
T_{\bbL} \XXm / T_{\bbL} \mathcal{Z}_{p,q} \simeq \Sym(\bbL^{\ast}) / \Hom(\bbL_V, V/\bbL_V)
\simeq \Sym(\bbL_V^{\ast}) \oplus \Sym(\bbL_{V^{\ast}}^{\ast})
\end{equation*}
Thus, the fiber of the conormal bundle $ T_{\mathcal{Z}_{p,q}}^{\ast} \XXm $ at $ \bbL $ is 
$ \Sym(\bbL_V) \oplus \Sym(\bbL_{V^{\ast}}) $.  

On the other hand, the fiber of $ \Tbundle_p $ at $ \bbL $ is $ \bbL_V $, while 
that of $ \Qbundle_p^{\ast} $ is $ (V/\bbL_V)^{\ast} \simeq \bbL_{V^{\ast}} $.  
Therefore the fiber of $ \Sym( \Qbundle_p^{\ast} ) \times \Sym( \Tbundle_p ) $ is naturally identified with 
$ \Sym(\bbL_V) \oplus \Sym(\bbL_{V^{\ast}}) $.
\end{proof}

For $ p + q = n $, let us consider the nilpotent $ \KGL $-orbit $ \orbit_{p,q} $ 
which is lifted from the trivial $ \MOO $-orbit $ \{ 0 \} $ in $ \lie{p}' $.  
Their shapes are all the same, namely it is $ 2^n $, and 
their signed Young diagrams are $ (+-)^p (-+)^q $.

\begin{theorem}
\label{theorem:resolution.of.nilpotent.orbit}
The affine quotient $ \nullcone_{p,q} \git \MOO $ is isomorphic to $ \closure{\orbit}_{p,q} $.  
Moreover, the image of the natural moment map $ \mu : T_{\mathcal{Z}_{p,q}}^{\ast} \XXm \to \lie{p} $ coincides with $ \closure{\orbit}_{p,q} $, 
and $ \mu $ gives a resolution of singularities of $ \closure{\orbit}_{p,q} $.  
Thus we have the following commutative diagram.
\begin{equation*}
\xymatrix @C+30pt {
& \ar[dl]_{\git \MOO} \makebox[20pt][c]{$\Hom( \Qbundle_p , \mathcal{U} ) \times \Hom( \mathcal{U} , \Tbundle_p )$} \ar[dr]^{\git K} \ar[d]^{\nu} & \\
\rule{0pt}{4ex}\makebox[10pt][c]{$T_{\mathcal{Z}_{p,q}}^{\ast} \XXm $} \ar[d]_{\mu} & 
\ar[dl]^{\git \MOO} {\nullcone}_{p,q} \ar[dr]_{\git K} & \{ \text{\upshape{pt}} \} \ar[d]\\
\closure{\orbit}_{p,q} &   & \{ 0 \} 
}
\end{equation*}
\end{theorem}

\begin{proof}
Let us denote $ X = \nullcone_{p,q} \git \MOO $.  
Take subspaces $ U^p \subset U $ and $ U^{\ast q} \subset U^{\ast} $ of dimension $ p $ and $ q $ respectively, 
and assume that the symmetric form restricted to these spaces are non-degenerate.  
If we put $ W^{p,q} = \Hom(V, U^p) \times \Hom(U^q , V ) $, 
then the image of $ \nullcone_{p,q} \cap W^{p,q} $ by the moment map $ \varphi $ is $ \closure{\orbit}_{p,q} $ 
by Theorem \ref{theorem:NOZ} (see also \S \ref{subsection:reduction.to.stable.range}).  
Thus $ X $ contains $ \closure{\orbit}_{p,q} $ as a closed subvariety.  
On the other hand, from Theorem \ref{theorem:decomposition.of.harmonics}, we can calculate
\begin{equation*}
\C[X] = \C[\nullcone_{p,q}]^{\MOO}
\simeq 
\bigoplus_{\alpha \in \partition_p, \; \beta \in \partition_q}
(\rho^{(n) \ast}_{\alpha} \boxtimes \rho^{(n)}_{\beta})^{\MOO} \otimes \rho^{(n)}_{\alpha \composit \beta} 
\simeq
\bigoplus_{\alpha \in \partition_p, \; \beta \in \partition_q}
\rho^{(n)}_{2 \alpha \composit 2 \beta} 
\end{equation*}
as a $ \KGL $-module.  
The same calculation leads the conclusion 
that it is isomorphic to $ \C[\closure{\orbit}_{p,q}] $ (see also \cite{Nishiyama.Tuebingen.2000}).  
Therefore we must have $ X = \closure{\orbit}_{p,q} $.

Let us prove that $ \mu $ is a resolution of $ \closure{\orbit}_{p,q} $.  
Since the conormal bundle is a bundle over a projective variety $ \mathcal{Z}_{p,q} $, $ \mu $ is proper.  
Take $ x \in \orbit_{p,q} $.  Then it is easy to see that the fiber $ \varphi^{-1}(x) $ is a single $ \MOO $-orbit, 
which is contained in $ \nullconeo_{p,q} $.  
Since $ \nu $ restricted to $ \nu^{-1}(\nullconeo_{p,q}) \to \nullconeo_{p,q} $ is isomorphism, 
the fiber $ \mu^{-1}(x) $ consists of a single point.  Thus $ \mu $ is birational.  
The closure $ \closure{\orbit}_{p,q} $ contains a nilpotent orbit $ \orbit_{s,t} \; ( s \leq p , t \leq q ) $, 
and no other orbit can adherent to $ \orbit_{p,q} $.  
For $ y \in \orbit_{s,t} $, it is easy to see that $ \mu^{-1}(y) $ has positive dimension.
\end{proof}

\section{Case of $ \OO^{\ast}_{2n} \times \Sp_{n,n} $}
\label{section:O2n.ast.Spnn}

Throughout in this section we assume that $ n $ is even.  
Then the results in \S \ref{section:Sp.O} have analogue for the pair $ \OO^{\ast}_{2n} \times \Sp_{n,n} $.  
We only state the results here because the proof is similar to the case of the pair $ ( \Sp_{2n}(\R), \OO_{n,n}(\R) ) $.

We consider a symplectic form on $ U (= U_1) $, and define the symplectic group $ \Sp(U) $ with respect to the form.  
By this non-degenerate symplectic form, we can identify $ U $ and $ U^{\ast} $, 
and we put $ U_2 = U^{\ast} $.  

We define $ \HGLxGL = \GL(U) \times \GL(U^{\ast}) $ as above, and put $ \MOO := \Sp(U) \times \Sp(U^{\ast}) \subset \HGLxGL $.  
Note that $ \MOO $ is the complexification of a maximal compact subgroup $ \USp_n \times \USp_n $ of 
the indefinite symplectic group $ \Sp_{n,n} $ over Hamiltonian numbers, where $ \USp_n $ denotes the unitary symplectic group of size $ n $ 
(the rank is $ n / 2 $).  
Similarly, in this section, 
$ \KGL = \GL_n(\C) $ should be regarded as the complexification of a maximal compact subgroup $ \U_n(\C) $ of $ \OO^{\ast}_{2n} $.
Thus we are working on the dual pair $ ( \OO^{\ast}_{2n} , \Sp_{n,n} ) $.

\begin{theorem}
Let us denote the conormal bundle of $ \mathcal{Z}_{p,q} \subset \XXp $ by $ T_{\mathcal{Z}_{p,q}}^{\ast} \XXp $.  
Then there exists a $ \KGL $-equivariant isomorphism 
\begin{equation*}
T_{\mathcal{Z}_{p,q}}^{\ast} \XXp \simeq \Alt( \Qbundle_p^{\ast} ) \times \Alt( \Tbundle_p ) \simeq \resolution_{p,q} \GITquotient \MOO  ,
\end{equation*}
where $ \Alt(X) $ denotes the alternating tensors over $ X $ of the second degree.  
Thus the conormal bundle of the closed $ \KGL $-orbit $ \mathcal{Z}_{p,q} $ in $ \XXp $ is obtained by 
taking the categorical quotient of the resolution of singularities $ \resolution_{p,q} $ of $ \nullcone_{p,q} $.
\end{theorem}

Let us consider nonnegative \emph{even} integers $ p $ and $ q $ such that $ p + q = n $.  
Then there is a nilpotent $ \KGL $-orbit $ \orbit_{p,q} $ 
which is lifted from the trivial $ \MOO $-orbit $ \{ 0 \} $ in $ \lie{p}' $.  
Their shapes are all the same, namely it is $ 2^n $, and 
their signed Young diagrams are $ (+-)^p (-+)^q $.

\begin{theorem}
Assume that $ p, q $ are even and satisfy $ p + q = n $. 
Then the affine quotient $ \nullcone_{p,q} \git \MOO $ is isomorphic to $ \closure{\orbit}_{p,q} $.  
Moreover, the image of the moment map $ \mu : T_{\mathcal{Z}_{p,q}}^{\ast} \XXp \to \lie{p} $ coincides with $ \closure{\orbit}_{p,q} $, 
and $ \mu $ gives a resolution of singularities of $ \closure{\orbit}_{p,q} $.  
Thus we have the following commutative diagram.
\begin{equation*}
\xymatrix @C+30pt {
& \ar[dl]_{\git \MOO} \makebox[20pt][c]{$\Hom( \Qbundle_p , \mathcal{U} ) \times \Hom( \mathcal{U} , \Tbundle_p )$} \ar[dr]^{\git K} \ar[d]^{\nu} & \\
\rule{0pt}{4ex}\makebox[10pt][c]{$T_{\mathcal{Z}_{p,q}}^{\ast} \XXp $} \ar[d]_{\mu} & 
\ar[dl]^{\git \MOO} {\nullcone}_{p,q} \ar[dr]_{\git K} & \{ \text{\upshape{pt}} \} \ar[d]\\
\closure{\orbit}_{p,q} &   & \{ 0 \} 
}
\end{equation*}
\end{theorem}

Let us take even integers $ p $ and $ q $ such that $ p + q = n $ as above.  
Then $ \nullcone_{p + 1, q - 1} \git \MOO = \closure{\orbit}_{p, q - 2} $ which is adherent to $ \orbit_{p,q} $.  
Such kind of degeneration of the categorical quotient makes the situation more complicated than 
the case of dual pair $ ( \Sp_{2 n}(\R) , \OO_{n, n}(\R) ) $.



\def\cftil#1{\ifmmode\setbox7\hbox{$\accent"5E#1$}\else
  \setbox7\hbox{\accent"5E#1}\penalty 10000\relax\fi\raise 1\ht7
  \hbox{\lower1.15ex\hbox to 1\wd7{\hss\accent"7E\hss}}\penalty 10000
  \hskip-1\wd7\penalty 10000\box7}
\providecommand{\bysame}{\leavevmode\hbox to3em{\hrulefill}\thinspace}

\end{document}